\documentclass[reqno]{amsart}

\usepackage{amssymb}

\newtheorem{theorem}{Theorem}[section]
\newtheorem{lemma}[theorem]{Lemma}

\newtheorem{corollary}[theorem]{Corollary}
\theoremstyle{definition}
\newtheorem{definition}[theorem]{Definition}

\theoremstyle{remark}
\newtheorem{remark}[theorem]{Remark}

\numberwithin{equation}{section}

\begin{document}

\title[Averaging in Differential Equations]
{Averaging Theorems for Ordinary Differential Equations
and Retarded Functional Differential Equations}

%    author one information
\author{Mustapha Lakrib}
\address{Laboratoire de Math\'ematiques, Universit\'e Djillali Liab\`es, B.P. 89, 22000 Sidi Bel Abb\`es, Alg\'erie}
\email{mlakrib@univ-sba.dz}

%    author two information
\author{Tewfik Sari}
\address{Laboratoire de Mathématiques, Informatique et Applications,
4, rue des Fr\`eres Lumi\`ere, 68093 Mulhouse, France}
\email{Tewfik.Sari@uha.fr}

%    \subjclass is required.

\date{}

%    Abstract is required.
\begin{abstract}
We prove averaging theorems for ordinary
differential equations
and retarded functional differential equations. Our assumptions are weaker 
than those required in the results of the existing literature. Usually,
we require that the nonautonomous
differential equation and the autonomous averaged equation
are locally Lipschitz and that 
the solutions of both equations exist on some interval. We extend
this result to the case of vector fields which are continuous in
the spatial variable uniformly with respect to time and without
any assumption on the interval of existence of the solutions of
the nonautonmous differential equation. 
Our results are formulated in classical
mathematics. Their proofs use nonstandard analysis.
\end{abstract}

\maketitle

\section{Introduction}

Averaging is an important method for analysis of
nonlinear oscillation equations containing a small parameter. This
method is well-known for ordinary differential equations (in short
ODEs) and fundamental averaging results (see, for instance,
\cite{Bogolyubov-Mitropolsky,Guckenheimer-Holmes,Hale2,Sanders-Verhulst,Sari1}
and references therein) assert that the solutions of a
nonautonomous equation in \emph{normal form}
\begin{equation}
\label{equation11}
 x'(\tau)=\varepsilon f(\tau,x(\tau)),
\end{equation}
where $\varepsilon$ is a small positive parameter, are
approximated by the solutions of the autonomous averaged equation
\begin{equation}
\label{equation12}
 y'(\tau)=\varepsilon F(y(\tau)).
\end{equation}
The approximation holds on time intervals of order $1/\varepsilon$
when $\varepsilon$ is sufficiently small. In (\ref{equation12}),
the function $F$ is the average of the function $f$ in
(\ref{equation11}) defined by
\begin{equation}\label{average}
F(x)=\lim_{T\rightarrow\infty}\frac{1}{T}\int_0^Tf(t,x)dt.
\end{equation}

The method of averaging was extended by Hale \cite{Hale1}  
(see
also Section 2.1 of \cite{Lehman}) 
to the
case of retarded functional differential equations (in short
RFDEs) containing a small parameter when the equations are
considered in normal form
\begin{equation}
\label{equation13}
 x'(\tau)=\varepsilon f(\tau,x_\tau),
\end{equation}
where, for $\theta\in[-r,0]$, $x_\tau(\theta)=x(\tau+\theta)$. 
Equations of the form (\ref{equation13})
cover a wide class of differential equations including those with
pointwise delay for which a method of averaging was developed in
\cite{Halanay,Medvedev,Volosov-Medvedev-Morgunov}. Note that the
averaged equation corresponding to (\ref{equation13}) is the ODE
\begin{equation}
\label{equation14}
 y'(\tau)=\varepsilon F(\tilde{y}^\tau),
\end{equation}
where, for $\tau$ fixed and $\theta\in[-r,0]$,
$\tilde{y}^\tau(\theta)=y(\tau)$
and the average function $F$ is defined by (\ref{average}). 
Recently, Lehman and Weibel
\cite{Lehman-Weibel} proposed to retain the delay in the averaged
equation and proved that equation
(\ref{equation13}) is approximated by the averaged RFDE
\begin{equation}\label{lwa}
 y'(\tau)=\varepsilon F(y_\tau).
\end{equation}
They observed, using numerical simulations, that equation
(\ref{equation13}) is better approximated by the averaged RFDE (\ref{lwa})
than by the averaged ODE (\ref{equation14}). However, this equation depends nontrivially
on the small parameter $\varepsilon$ (see Remark~\ref{rem3}).

The change from the slow time scale $\tau$ to the fast time scale
$t=\tau/\varepsilon$ transforms equations (\ref{equation11}) and
(\ref{equation12}), respectively, into
\begin{equation}
\label{equation15}
 \dot{x}(t)=f(t/\varepsilon ,x(t))
\end{equation}
and
\begin{equation}
\label{equation16}
 \dot{y}(t)=F(y(t)).
\end{equation}
Thus a method of averaging can be developed for
(\ref{equation15}), that is, if $\varepsilon$ is sufficiently
small, the difference between the solution $x$ of
(\ref{equation15}) and the solution $y$ of (\ref{equation16}),
with the same initial condition, is small on finite time
intervals.

The analog of equation (\ref{equation15}) for RFDEs is
\begin{equation}
\label{equation17}
 \dot{x}(t)=f\left(t/\varepsilon,x_t\right).
\end{equation}
The averaged equation corresponding to (\ref{equation17}) is the
RFDE
\begin{equation*}
 \dot y(t)=F(y_t),
\end{equation*}
where the average function $F$ is defined by (\ref{average}).

Notice that the RFDEs (\ref{equation13}) and (\ref{equation17}) are not 
equivalent under the change of time $t=\tau/\varepsilon$, as it was the case 
for the ODEs (\ref{equation11}) and (\ref{equation15}). Indeed, by rescaling
$\tau$ as  $t=\tau/\varepsilon$ equation (\ref{equation13}) becomes
\begin{equation}
\label{equation18}
 \dot{x}(t)=\displaystyle f\left(t/\varepsilon,x_{t,\varepsilon}\right),
\end{equation}
where, for $\theta\in[-r,0]$,
$x_{t,\varepsilon}(\theta)=x(t+\varepsilon\theta)$. Equation
(\ref{equation18}) is different from~(\ref{equation17}), so that
the results obtained for (\ref{equation18}) cannot be applied to
(\ref{equation17}). This last equation deserves a special
attention. It was considered by Hale and Verduyn Lunel  in
\cite{Hale-Lunel1} where a method of averaging is developed for
infinite dimensional evolutionary equations which include RFDEs
such (\ref{equation17}) as a particular case (see also Section
12.8 of Hale and Verduyn Lunel's book \cite{Hale-Lunel2} and
Section~2.3 of \cite{Lehman}). 

Following our previous works
\cite{Lakrib1,Lakrib2,LakibThesis,Lakrib-Sari,Sari1,Sari2}, we
consider in this paper all equations (\ref{equation15}),
(\ref{equation17}) and (\ref{equation18}). Our aim is to give
theorems of averaging under weaker conditions  than those of the
literature.
We want to emphasize that our main contribution is the weakening of the 
regularity conditions on the equation under which the averaging method is justified in the existing
literature. Indeed, usually classical averaging
theorems require that the vector field $f$ in (\ref{equation15}),
(\ref{equation17}) and (\ref{equation18})  is at least locally
Lipschitz with respect to the second  variable uniformly with
respect to the first one (see Remarks \ref{rem1}, \ref{rem2} and \ref{rem4} below). In our results this condition is
weakened  and it is only assumed that $f$ is continuous in the second variable uniformly
with respect to the first one. Also, it is often assumed that the
solutions $x$ and $y$ exist on the same finite interval of time. In
this paper we assume only that the solution $y$ of the averaged
equation exists on some finite interval and we give conditions on
the vector field $f$ so that, for $\varepsilon$ sufficiently
small, the solution $x$ of (\ref{equation15}), (\ref{equation17}) or
(\ref{equation18}) will be defined at least on the same interval.
The {\em uniform quasi-boundedness} of the vector field $f$ is
thus introduced for this purpose. Recall that the property of
quasi-boundedness is strongly related to results on continuation
of solutions of RFDEs. It should be noticed that the existing literature 
\cite{Hale1,Hale-Lunel1,Lehman} proposed also important results on the infinite time interval
$[0,\infty)$, provided that more hypothesis are made on the nonautonomous system and its averaged system.
For example, to a hyperbolic equilibrium point of the averaged system there corresponds a 
periodic solution of the original equation if $\varepsilon$ is small. Of course, for such results, 
stronger assumptions on the regularity of the vector field $f$ are required.

In this work our averaging results are formulated in classical
mathematics. We prove them within \emph{Internal Set
Theory} (in short IST) \cite{Nelson} which is an axiomatic
approach to \emph{Nonstandard Analysis} (in short NSA)
\cite{Robinson}. The idea to use NSA in perturbation theory of
differential equations goes back to the 1970s with the Reebian
school \cite{Lutz,Lutz-Sari}. It has become today a
well-established tool in asymptotic theory, as attested by the
the five-digits classification 34E18 of the 2000 Mathematical
Subject Classification (see also
\cite{vdb,DD,Lakrib-Sari,Lobry-Sari-Touhami,Sari-Yadi}).

The structure of the paper is as follows. In Section
\ref{Notations and Main Results} we introduce the notations and
present our main results  : Theorems \ref{theorem21}, \ref{theorem24}
and~\ref{theorem28}. We discuss also both periodic and almost periodic
special cases. In Section~\ref{The Stroboscopic Method} we start with a short 
tutorial to NSA and then state our 
main (nonstandard) tool, the so-called
\emph{stroboscopic method}. 
In Section \ref{Proofs of the Results}.
we give the proofs of
Theorems~\ref{theorem21}, \ref{theorem24} and~\ref{theorem28}.

We recall that the stroboscopic method was proposed for the first
time in the study of some ODEs with a small parameter which occur
in the theory of nonlinear oscillations
\cite{Callot-Sari,Lutz,Sari2}. Here, we first present a slightly
modified version of this method  and then extend it in the context
of RFDEs.

Let us notice that none of our
proofs needs to be translated into classical mathematics, because IST is a
conservative extension of ordinary mathematics, that is, any
classical statement which is a theorem of IST is also a
theorem of ordinary mathematics.

\section{Notations and Main Results}
\label{Notations and Main Results}

In this section we will present our main results on averaging for
fast oscillating ODEs (\ref{equation15}), 
RFDEs in normal form (\ref{equation18}) and 
fast oscillating RFDEs (\ref{equation17}). First we introduce some
necessary notations.
We assume that $r\geq 0$ is a fixed real
number and denote by
$\mathcal{C}=\mathcal{C}([-r,0],\mathbb{R}^d)$ the Banach space of
continuous functions from $[-r,0]$ into $\mathbb{R}^d$ with the
norm $\|\phi\|=\sup\{|\phi(\theta)|: \theta\in[-r,0]\}$, where
$|\cdot|$ is a norm of $\mathbb{R}^d$. Let $L\geq 0$. If
$x:[-r,L]\rightarrow\mathbb{R}^d$ is a continuous function then,
for each $t\in[0,L]$, we define $x_t\in \mathcal{C}$ by setting
$x_t(\theta)=x(t+\theta)$ for all $\theta\in [-r,0]$. Note that
when $r=0$ the Banach space $\mathcal{C}$ can be identified with
$\mathbb{R}^d$ and $x_t$ with $x(t)$ for each $t\in[0,L]$.

\subsection{Averaging for ODEs}
\label{Averaging for ODEs}

Let $f:\mathbb{R}_+\times\mathbb{R}^d\rightarrow\mathbb{R}^d$,
$(t,x)\mapsto f(t,x)$, be a continuous function. Let
$x_0\in\mathbb{R}^d$ be an initial condition. We consider the initial
value problem
\begin{equation}
\label{equation21}
 \dot{x}(t)=f\left(t/\varepsilon,x(t)\right),\qquad x(0)=x_0,
\end{equation}
where $\varepsilon>0$ is a small parameter.
We state the precise assumptions on this problem.

\begin{enumerate}
    \item[(C1)]
The function $f$ is continuous in the second variable uniformly
with respect to the first one.
    \item[(C2)]
For all $x\in\mathbb{R}^d$, there exists a limit
$\displaystyle
F(x):=\lim_{T\rightarrow\infty}\frac{1}{T}\int_0^Tf(t,x)dt$.
\end{enumerate}

From conditions (C1) and (C2) we deduce that  the
average of the function $f$, that is, the function
$F:\mathbb{R}^d\rightarrow\mathbb{R}^d$ in \mbox{(C2)}, is
continuous (see Lemma \ref{lemma41}). So, the following (averaged)
initial value problem is well defined.
\begin{equation}
\label{equation22}
 \dot y(t)=F(y(t)),\qquad y(0)=x_0.
\end{equation}

We need also the condition :
\begin{enumerate}
    \item[(C3)]
The initial value problem (\ref{equation22}) has a unique
solution.
\end{enumerate}

The main theorem of this section is on averaging for fast
oscillating ODEs. It
establishes nearness of solutions of (\ref{equation21}) and
(\ref{equation22}) on finite time intervals, and reads as follows.

%%%%%%%%%%%%%%%
\begin{theorem}
\label{theorem21} Let
$f:\mathbb{R}_+\times\mathbb{R}^d\rightarrow\mathbb{R}^d$ be a
continuous function and let $x_0\in\mathbb{R}^d$. Suppose  that
conditions (C1)-(C3) are satisfied. Let $y$ be the solution of
(\ref{equation22}) and let $L\in J$, where $J$ is the positive
interval of definition of $y$. Then, for every $\delta>0$, there
exists $\varepsilon_0=\varepsilon_0(L,\delta)>0$ such that, for
all $\varepsilon\in(0,\varepsilon_0]$, every solution $x$ of
(\ref{equation21}) is defined at least on $[0,L]$ and satisfies
$|x(t)-y(t)|<\delta$ for all $t\in [0,L]$.
\end{theorem}
%%%%%%%%%%%%%

Let us discuss now the result above when  the function $f$ is
periodic or more generally almost periodic in the first variable.
We will see that some of the conditions in Theorem \ref{theorem21}
can be removed.
Indeed, in the case where $f$ is periodic in $t$, from continuity plus
periodicity properties one can easily deduce
condition (C1). Periodicity also implies  condition (C2) in an
obvious way. The average of $f$ is then given, for every
$x\in\mathbb{R}^d$, by 
\begin{equation}\label{periodic}
F(x)=\frac{1}{T}\int_0^Tf(t,x)dt,
\end{equation}
where $T$ is the period.
In the case where $f$ is almost periodic in $t$ it is well-known that for 
all $x\in\mathbb{R}^d$, the limit
\begin{equation}\label{almostperiodic}
F(x)=\lim_{T\rightarrow\infty}\frac{1}{T}\int_s^{s+T}f(t,x)dt
\end{equation}
exists uniformly with respect to $s\in\mathbb{R}$. 
So, condition (C2) is satisfied when $s=0$.
We point out also that in a number of cases encountered in applications
the function $f$ is a finite sum of periodic functions in $t$. As in the
periodic case above, condition~(C1) is satisfied. Hence we have the following
result.

%%%%%%%%%%%%%%%
\begin{corollary}[Periodic and Almost periodic cases]
The conclusion of Theorem \ref{theorem21} holds when 
$f:\mathbb{R}_+\times\mathbb{R}^d\rightarrow\mathbb{R}^d$ is a
continuous function 
which is periodic (or a sum of periodic functions) 
in the first variable and satisfies condition (C3).
It holds also when $f$ is continuous,
almost periodic in the first variable and satisfies conditions
(C1) and (C3).
\end{corollary}
%%%%%%%%%%%%%

\begin{remark}\label{rem1}
In the results of the classical literature, 
for instance \cite[Theorem 1, p. 202]{Lehman}, 
it is assumed that $f$ has bounded partial
derivatives with respect to the second variable.
\end{remark}

\subsection{Averaging for RFDEs in normal form}
\label{Averaging for RFDEs in normal form}

This section concerns the use of the method of averaging to
approximate initial value problems of the form
\begin{equation}
\label{equation23}
 \dot{x}(t)= f(t/\varepsilon,x_{t,\varepsilon}),
\quad x(t)=\phi(t/\varepsilon), t\in[-\varepsilon r,0].
\end{equation}
Here $f:\mathbb{R}_+\times\mathcal{C}\rightarrow\mathbb{R}^d$,
$(t,x)\mapsto f(t,x)$, is a continuous function, $\phi\in
\mathcal{C}$ is an initial condition and $\varepsilon>0$ is a
small parameter. For each $t\geq 0$, $x_{t,\varepsilon}$ denotes
the element of $\mathcal{C}$ given by
$x_{t,\varepsilon}(\theta)=x(t+\varepsilon\theta)$ for all
$\theta\in[-r,0]$.

We recall that 
the change of time scale $t=\varepsilon\tau$ transforms 
(\ref{equation23}) into the following initial
value problem, associated to a RFDE in normal form :
\begin{equation*}
x'(\tau)=\displaystyle\varepsilon f(\tau,x_\tau), \quad x_0=\phi.
\end{equation*}

We make the following hypotheses:

\begin{enumerate}
    \item[(H1)]
The function $f$ is continuous in the second variable uniformly
with respect to the first one.
    \item[(H2)]
The function $f$ is  quasi-bounded in the second variable
uniformly with respect to the first one, that is, for every
bounded subset $B$ of $\mathcal{C}$, $f$ is bounded on
$\mathbb{R}_+\times B$.
    \item[(H3)]
For all $x\in \mathcal{C}$, the limit $\displaystyle
F(x):=\lim_{T\rightarrow\infty}\frac{1}{T}\int_0^Tf(t,x)dt$
exists.
\end{enumerate}

We define the averaged initial value problem associated to
(\ref{equation23}) by
\begin{equation}
\label{equation24}
 \dot y(t)=G(y(t)),\qquad y(0)=\phi(0).
\end{equation}
The function $G:\mathbb{R}^d\rightarrow\mathbb{R}^d$ is defined by
$G(x)=F(\tilde{x})$ where, for each
$x\in\mathbb{R}^d$, $\tilde{x}\in\mathcal{C}$ is given by
$\tilde{x}(\theta)=x$, $\theta\in[-r,0]$.
We add the last hypothesis:
\begin{enumerate}
    \item[(H4)]
The initial value problem (\ref{equation24}) has a unique
solution.
\end{enumerate}

As we will see later, condition (H2) is used essentially to prove
continuability of solutions of (\ref{equation23}) at least on
every finite interval of time on which the solution of
(\ref{equation24}) is defined. For more details and a complete
discussion about quasi-boundedness and its crucial role in the
continuability of solutions of RFDEs, we refer the reader to
Sections 2.3 and 3.1 of \cite{Hale-Lunel2}.

In assumption (H4) we anticipate the existence of solutions of
(\ref{equation24}). This will be justified a posteriori by Lemma
\ref{lemma41} where we show that the function
$F:\mathcal{C}\rightarrow\mathbb{R}^d$ in (H3), which is the
average of the function $f$, is continuous. This implies the
continuity of $G:\mathbb{R}^d\rightarrow\mathbb{R}^d$ in
(\ref{equation24}) and then guaranties the existence of solutions.

The result below is our main theorem on averaging for RFDEs in
normal form. It states closeness of solutions of
(\ref{equation23}) and (\ref{equation24}) on finite time
intervals.

%%%%%%%%%%%%%%%
\begin{theorem}
\label{theorem24} Let
$f:\mathbb{R}_+\times\mathcal{C}\rightarrow\mathbb{R}^d$ be a
continuous function and $\phi\in\mathcal{C}$. Let conditions
(H1)-(H4) hold. Let $y$ be the solution of (\ref{equation24}) and
let  $L\in J$, where $J$ is the positive interval of definition of
$y$. Then, for every $\delta>0$, there exists
$\varepsilon_0=\varepsilon_0(L,\delta)>0$ such that, for all
$\varepsilon\in(0,\varepsilon_0]$, every solution $x$ of
(\ref{equation23}) is defined at least on $[-\varepsilon r,L]$ and
satisfies $|x(t)-y(t)|<\delta$ for all $t\in [0,L]$.
\end{theorem}
%%%%%%%%%%%%%

As in Section \ref{Averaging for ODEs}, we discuss now both periodic  
and almost periodic special cases. In
each one, some of the conditions in Theorem~\ref{theorem24} can be
either removed or weakened. Let us consider the following (weak) condition
which will be used hereafter instead of condition (H2):

\begin{enumerate}
    \item[(H5)]
The function $f$ is  quasi-bounded, that is, $f$ is bounded on
bounded subsets of $\mathbb{R}_+\times\mathcal{C}$.
\end{enumerate}

When $f$ is periodic it is easy to see that  condition (H1) derives from the continuity
and the periodicity properties of $f$. On the other
hand, by periodicity and condition (H5), condition
(H2) is also satisfied. The average $F$ in condition (H3) exists and is now given by formula
(\ref{periodic}) where $T$ is the period.
When $f$ is almost periodic, condition (H5) imply condition (H2) and the average $F$ is given by formula
(\ref{almostperiodic}).
Quite often the function $f$ is a finite sum of periodic functions
so that condition (H1) is satisfied. Hence we have the following
result.

%%%%%%%%%%%%%%%
\begin{corollary}[Periodic and Almost periodic cases]
\label{theorem25}
The conclusion of Theorem \ref{theorem24} holds when 
$f:\mathbb{R}_+\times\mathcal{C}\rightarrow\mathbb{R}^d$
is a continuous function 
which is periodic (or a sum of periodic functions) 
in the first variable and satisfies condition (H4) and (H5).
It holds also when $f$ is continuous,
almost periodic in the first variable and satisfies conditions
(H1), (H4) and (H5).
\end{corollary}

Consider now the special case of equations with pointwise delay of the form
\begin{equation*}
\dot{x}(t)=f(t/\varepsilon,x(t), x(t-\varepsilon r))
\end{equation*}
which is obtained, by letting $\tau=t/\varepsilon$, from equation
\begin{equation*}
x'(\tau)=\varepsilon f(\tau,x(\tau), x(\tau-r)).
\end{equation*}
In this case, 
for both periodic and almost
periodic functions, condition
(H5) follows from the continuity property and then may
be removed in Corollary \ref{theorem25}.

\begin{remark}\label{rem2}
In the results of the
literature, for instance \cite[Theorem 3, p. 206]{Lehman}, 
$f$ is assumed to be locally Lipschitz with
respect to the second variable.
Note that local Lipschitz condition with respect to the second
variable implies condition (H1). It also assures the local
existence for the solution of (\ref{equation23}). But, in
opposition to the case of ODEs, it is well known (see Sections 2.3
and 3.1 of \cite{Hale-Lunel2}) that without condition
(H5) one cannot extend the solution $x$ to finite time
intervals where the solution $y$ is defined in spite of the
closeness of $x$ and $y$. So, in the existing literature it is 
assumed that the solutions $x$ and $y$ are both defined
at least on the same interval $[0,L]$. 
\end{remark}

\begin{remark}\label{rem3}
In the introduction, we noticed that
Lehman and Weibel \cite{Lehman-Weibel} proposed to retain the delay in the averaged
equation (\ref{lwa}). At time scale
$t=\varepsilon\tau$, their observation is that equation
(\ref{equation23}) is better approximated by the averaged RFDE 
\begin{equation}\label{lwat}
 \dot{y}(t)=F(y_{t,\varepsilon})
\end{equation}
than by the averaged ODE (\ref{equation24}). It should be noticed that 
the averaged RFDE (\ref{lwat}) depends on the small parameter $\varepsilon$, which is not 
the case of the averaged equation (\ref{equation24}).
\end{remark}

\subsection{Averaging for fast oscillating RFDEs}

The aim here is to approximate the solutions of the initial value
problem
\begin{equation}
\label{equation25}
 \dot{x}(t)=\displaystyle
f\left(t/\varepsilon,x_t\right),\qquad x_0=\phi,
\end{equation}
where $f:\mathbb{R}_+\times\mathcal{C}\rightarrow\mathbb{R}^d$,
$(t,x)\mapsto f(t,x)$, is a continuous function, $\phi\in
\mathcal{C}$ is an initial condition and $\varepsilon>0$ is a
small parameter. This will be obtained under conditions (H1), (H2)
and (H3) in Section \ref{Averaging for RFDEs in normal form}
plus condition (H6) below. We define the averaged initial value problem associated to 
(\ref{equation25}) by
\begin{equation}
\label{equation26}
 \dot y(t)=F(y_t),\qquad y_0=\phi.
\end{equation}
As in the previous section, conditions (H1) and (H3) imply the
continuity of the function $F:\mathcal{C}\rightarrow\mathbb{R}^d$
in (H3). So, the problem (\ref{equation26}) is well defined.
We need the following condition
\begin{enumerate}
    \item[(H6)]
The initial value problem (\ref{equation26})
has a unique solution.
\end{enumerate}

Under the above assumptions, we may state our main result on
averaging for fast oscillating RFDEs.
It shows that the solution of (\ref{equation26}) is an
approximation of solutions of (\ref{equation25})
on finite time intervals.

%%%%%%%%%%%%%%%
\begin{theorem}
\label{theorem28} Let
$f:\mathbb{R}_+\times\mathcal{C}\rightarrow\mathbb{R}^d$ be a
continuous function and let $\phi\in\mathcal{C}$. Suppose that
conditions (H1)-(H3) and (H6) are satisfied. Let $y$
be the solution of (\ref{equation26}) and let  $L\in J$ be
positive, where $J$ is the interval of definition of $y$. Then,
for every $\delta>0$, there exists
$\varepsilon_0=\varepsilon_0(L,\delta)>0$ such that, for all
$\varepsilon\in(0,\varepsilon_0]$, every solution $x$ of
(\ref{equation25}) is defined at least on $[-r,L]$ and satisfies
$|x(t)-y(t)|<\delta$ for all $t\in [0,L]$.
\end{theorem}
%%%%%%%%%%%%%

In the same manner as in Section \ref{Averaging for RFDEs in
normal form} we have the following
result corresponding to the periodic and almost periodic special
cases.

\begin{corollary}[Periodic and Almost periodic cases]
\label{theorem29}
The conclusion of Theorem \ref{theorem28} holds when 
$f:\mathbb{R}_+\times\mathcal{C}\rightarrow\mathbb{R}^d$
is a continuous function 
which is periodic (or a sum of periodic functions) 
in the first variable and satisfies condition (H5) and (H6).
It holds also when $f$ is continuous,
almost periodic in the first variable and satisfies conditions
(H1), (H5) and (H6).
\end{corollary}

For fast oscillating equations with
pointwise delay of the form
\begin{equation*}
\dot{x}(t)=f(t/\varepsilon,x(t), x(t-r)),
\end{equation*}
in the periodic case as well as in the almost
periodic one, condition (H5) derives from the
continuity property and then can  be removed in
Corollary \ref{theorem29}.  

\begin{remark}\label{rem4}
In the results of the classical literature, for instance, \cite[Theorem 4, p. 210]{Lehman},
it is assumed that $f$ is locally Lipschitz with
respect to the second variable and the existence of the solutions 
$x$ and $y$ on the same interval $[0,L]$ is required.
\end{remark}

\section{The Stroboscopic Method}
\label{The Stroboscopic Method}

\subsection{Internal Set Theory}
\label{IST:A short tutorial}
In this section we give a short tutorial of NSA. 
Additional informations can be found in
\cite{vdb,DD,Nelson,Robinson}. 
\emph{Internal Set Theory} (IST) is a theory extending ordinary
mathematics, say ZFC (Zermelo-Fraenkel set theory with the axiom
of choice), that axiomatizes (Robinson's) \emph{nonstandard
analysis} (NSA). We adjoin a new undefined unary predicate
\emph{standard} (st) to ZFC. In addition to the usual axioms of
ZFC, we introduce three others for handling the new predicate in a
relatively \emph{consistent} way. Hence \emph{all theorems} of ZFC
\emph{remain valid} in IST. What is new in IST is an addition, not
a change.

A real number $x$ is  said to be \emph{infinitesimal} if $|x|<a$
for all standard positive real numbers $a$ and \emph{limited} if
$|x|\leq a$ for some standard positive real number $a$. A limited
real number which is not infinitesimal is said to be
\emph{appreciable}. A real number which  is not limited is said to
be  \emph{unlimited}. The  notations $x\simeq 0$ and
$x\simeq+\infty$ are used to denote, respectively, $x$ is
infinitesimal and  $x$ is unlimited positive.

Let $E$ be a standard normed space. 
A vector $x$ in $E$ is \emph{infinitesimal}
(resp. \emph{limited, unlimited}) if its norm $\|x\|$ is
infinitesimal (resp. limited, unlimited).
Two elements $x$ and $y$ in
$E$ are said to be \emph{infinitely close}, in symbols,  $x\simeq
y$, if $\|x-y\|\simeq 0$. 
An element $x$ is said to be \emph{nearstandard} if
$x\simeq x_0$ for some standard $x_0\in E$. The element $x_0$ is
called the \emph{standard part} or \emph{shadow} of $x$. It is
unique and is usually denoted by $^ox$.
Note that for $d$ standard,  each limited vector in $\mathbb{R}^d$ is
nearstandard.

Let $I\subset \mathbb{R}$ be some interval and
$f:I\rightarrow\mathbb{R}^d$ be a function, with $d$ standard. We
say that $f$ is \emph{S-continuous at} $x$ in $I$ if, for all $y$
in $I$, $x\simeq y$ implies $f(x)\simeq f(y)$, and
\emph{S-continuous on} $I$ if, for all  $x$ and $y$ nearstandard
in $I$, $x\simeq y$ implies $f(x)\simeq f(y)$. When $f$ (and then
$I$) and $x$ are standard, the first definition is the same as
saying that $f$ is continuous at $x$, and the uniform continuity
of $f$ on $I$ is equivalent to, for all $x$ and $y$ in $I$,
$x\simeq y$ implies $f(x)\simeq f(y)$.

We need the following basic result on S-continuous functions.

\begin{theorem}[\textmd{Continuous shadow}]
\label{theorem31} Let $I\subset \mathbb{R}$ be some interval and
$f:I\rightarrow\mathbb{R}^d$ be a function, with  $d$ a standard
positive integer. Let ${}^sI$ be the standard subset of
$\mathbb{R}$ whose standard elements are those of $I$. There
exists a standard and continuous function
$f_0:{}^sI\rightarrow\mathbb{R}^d$ such that, for all $x$
nearstandard in $I$, $f(x)\simeq f_0(x)$ holds, if and only if,
$f$ is S-continuous on $I$ and $f(x)$ is nearstandard for all $x$
nearstandard in $I$.
\end{theorem}

The function $f_0$ in Theorem \ref{theorem31} is unique. It is
called the \emph{standard part {\rm or} shadow} of the function
$f$ and denoted by ${}^of$.

\subsection{The Stroboscopic Method for ODEs}

Let $x_0\in\mathbb{R}^d$ be standard and let $F:\mathbb{R}_+\times
\mathbb{R}^d\rightarrow\mathbb{R}^d$ be a standard and continuous
function. Let $I$ be some subset of $\mathbb{R}$ and let
$x:I\rightarrow\mathbb{R}^d$ be a function such that $0\in I$ and
$x(0)\simeq x_0$.

%%%%%%%%%%%%%%%%%%
\begin{definition}[\emph{F}-Stroboscopic property]
\label{definition1} The function $x$ is said to satisfy the
\emph{F}-Stroboscopic property if there exists $\mu>0$ such that,
for all limited $t\geq 0$ in $I$ satisfying $[0,t]\subset I$ and
$x(s)$ is limited for all $s\in[0,t]$, there exists $t'\in I$ such
that $\mu<t'-t\simeq 0$, $[t,t']\subset I$,\ $x(s)\simeq x(t)$ for
all $s\in [t,t']$ and $\displaystyle \frac{x(t')-x(t)}{t'-t}\simeq
F(t,x(t)). $
\end{definition}
%%%%%%%%%%%%%%%%

Now, if a function satisfies the \emph{F}-stroboscopic property,
the result below asserts that it can be approximated by a solution
of the ODE
\begin{equation}
\label{equation31}
 \dot{y}(t) =F(t,y(t)),\qquad y(0)=x_0.
\end{equation}

%%%%%%%%%%%%%%%
\begin{theorem}[Stroboscopic Lemma for ODEs]
\label{theorem32} Suppose that\\
{\rm (a)} The function $x$ satisfies the F-stroboscopic property.\\
{\rm (b)} The initial value problem (\ref{equation31}) has a unique solution $y$.
Let $J=[0,\omega)$, $0<\omega\leq\infty$, be its  maximal positive
interval of definition.\\
Then, for every standard $L\in J$, $[0,L]\subset I$ and the
approximation \mbox{$x(t)\simeq y(t)$} holds for all $t\in [0,L]$.
\end{theorem}
%%%%%%%%%%%%%

The proof of Stroboscopic  Lemma for ODEs needs some results which
are given in the section below.

\subsubsection{Preliminaries}

%%%%%%%%%%%%%
\begin{lemma}
\label{lemma31} Let $L>0$ be limited such that $[0,L]\subset I$.
Suppose that\\
{\rm (i)} The function $x$ is limited on $[0,L]$.\\
{\rm (ii)}
There exist some positive integer $N$ and some infinitesimal
partition $\{t_n:n=0,\ldots,N+1\}$ of $[0,L]$ such that $t_0=0$,
$t_N\leq L<t_{N+1}$ and, for $n=0,\ldots,N$, $t_{n+1}\simeq t_n$,
$x(t)\simeq x(t_n)$ for all $t\in [t_n,t_{n+1}]$, and
$\displaystyle \frac{x(t_{n+1})-x(t_n)}{t_{n+1}-t_n}\simeq
F(t_n,x(t_n))$.\\
Then the function $x$ is S-continuous at each point in $[0,L]$.
\end{lemma}
%%%%%%%%%%%
\begin{proof}
Let $t\in[0,L]$. We will show that $x$ is S-continuous at $t$. Let
$t'\in[0,L]$ and $p,q\in\{0,\ldots,N\}$ be such that $t\leq t'$,
$t\simeq t'$, $t\in[t_p,t_{p+1}]$ and $t'\in[t_q,t_{q+1}]$. We
write
\begin{equation}
\label{equation32}
x(t_q)-x(t_p)  = 
\displaystyle\sum_{n=p}^{q-1}(x(t_{n+1})-x(t_n))
    = \displaystyle \sum_{n=p}^{q-1}
(t_{n+1}-t_n)[F(t_n,x(t_n))+\eta_n],
\end{equation}
where $\eta_n\simeq 0$ for all $n\in\{p,\ldots,q-1\}$. Denote
\begin{equation*}
\eta=\max_{p\leq n\leq q-1}|\eta_n|,\qquad
 m=\max_{p\leq n\leq q-1}|F(t_n,x(t_n))|.
\end{equation*}
We have $\eta\simeq 0$ and $m=|F(t_s,x(t_s))|$ for some
$s\in\{p,\dots,q-1\}$. Since $(t_s,x(t_s))$ is 
limited, it is nearstandard. 
Since the function $F$ is standard and
continuous, $F(t_s,x(t_s))$ is nearstandard. So is $m$.
Hence (\ref{equation32}) leads to the approximation
\begin{equation*}
|x(t')-x(t)|\simeq |x(t_q)-x(t_p)| \leq (m+\eta)(t_q-t_p) \simeq 0
\end{equation*}
 which proves the S-continuity of $x$ at $t$ and completes the
proof.
\end{proof}
%%%%%%%%%%%

When instead of $L$ limited we suppose $L$ standard, Lemma
\ref{lemma31} transforms to the following result with more
properties about the function $x$.

%%%%%%%%%%%%%
\begin{lemma}
\label{lemma32}
Let $L>0$ be standard such that $[0,L]\subset I$.
Suppose that conditions (i) and (ii) in Lemma \ref{lemma31} hold.
Then the standard function $y:[0,L]\rightarrow\mathbb{R}^d$
defined, for all standard $t\in[0,L]$, by
$y(t)={}^o\big(x(t)\big)$, is a solution of (\ref{equation31}).
Moreover, the approximation $x(t)\simeq y(t)$ holds for all
$t\in[0,L]$.
\end{lemma}
%%%%%%%%%%%
\begin{proof} To prove the lemma  we proceed in two steps.

    \emph{Step 1.}
We claim that the function $y$ is continuous on $[0,L]$. Indeed,
by Lemma \ref{lemma31} the function $x$ is S-continuous on
$[0,L]$. Taking hypothesis (i) into account, the claim follows
from Theorem \ref{theorem31}. We have moreover
\begin{equation*}
    y(t)\simeq y({}^ot)\simeq x({}^ot)\simeq x(t),\quad \forall t\in[0,L].
\end{equation*}
The first part of the proof is complete.

    \emph{Step 2.} To show that the function  $y$ satisfies, for all $t\in[0,L]$,
 \begin{equation*}
y(t)=x_0+\int_0^t F(s,y(s))ds,
\end{equation*}
let  $t\in [0,L]$ and $n\in\{0,\ldots,N\}$ be such that
$t\in[t_n,t_{n+1}]$ with $t$ standard. Then
\begin{eqnarray}
\label{equation33}
y(t)-x_0 \simeq  x(t_n)-x(0) &=&  \displaystyle
\sum_{k=0}^{n-1}(x(t_{k+1})-x(t_k))\nonumber\\
      & = & \displaystyle\sum_{k=0}^{n-1}(t_{k+1}-t_k)[F(t_k,x(t_k))+\eta_k],
\end{eqnarray}
where $\eta_k\simeq 0$ for all $k\in\{0,\ldots,n-1\}$. As $F$ is
standard and continuous, and by Step 1 above $x(t_k)\simeq y(t_k)$
with $x(t_k)$ nearstandard, we have $F(t_k,x(t_k))\simeq
F(t_k,y(t_k))$ so that (\ref{equation33}) gives
 \begin{equation*}
y(t)-x_0  \simeq  \displaystyle
\sum_{k=0}^{n-1}(t_{k+1}-t_k)[F(t_k,y(t_k))+\beta_k+\eta_k]
\simeq \displaystyle \int_0^t F(s,y(s))ds,
\end{equation*}
where $\beta_k\simeq 0$ for all $k\in\{0,\ldots,n-1\}$. Thus the
approximation
\begin{equation}
\label{equation34}
 y(t)\simeq x_0+ \int_0^t F(s,y(s))ds
\end{equation}
holds for all standard $t\in[0,L]$. Actually (\ref{equation34}) is
an equality since  both sides of which are standard. We have thus,
for all standard $t\in[0,L]$,
\begin{equation}
\label{equation35}
 y(t)=x_0+\int_0^t F(s,y(s))ds
\end{equation}
and by transfer (\ref{equation35}) holds for all $t\in[0,L]$. The
proof is complete.
\end{proof}
%%%%%%%%%%%

The following statement is a consequence of Lemma \ref{lemma32}.

%%%%%%%%%%%%%
\begin{lemma}
\label{lemma33}
Let $L>0$ be standard  such that $[0,L]\subset I$.
Suppose that\\
{\rm(i)}
The function $x$ is limited on $[0,L]$.\\
{\rm(ii)}
There exists $\mu>0$ such that, for all $t\in [0,L]$, there exists
$t'\in I$ such that $\mu<t'-t\simeq 0$, $[t,t']\subset I$,
$x(s)\simeq x(t)$ for all $s\in [t,t']$, and 
$\displaystyle\frac{x(t')-x(t)}{t'-t}\simeq F(t,x(t))$.\\
Then the  function $x$ is S-continuous on $[0,L]$ and its shadow
is a solution $y$ of (\ref{equation31}). So, we have $x(t)\simeq
y(t)$ for all $t\in[0,L]$
\end{lemma}
%%%%%%%%%%%
\begin{proof}
First of all we have $\lambda\in A_\mu$ for all standard real
number $\lambda>0$, where  $A_\mu$ is the subset of $\mathbb{R}$
defined by
\begin{equation*}
 A_\mu=\{\lambda\in\mathbb{R}\ /\
\forall t\in[0,L]\ \exists t'\in I:
\mathcal{P}_\mu(t,t',\lambda)\}
\end{equation*}
 and $\mathcal{P}_\mu(t,t',\lambda)$ is the property
$$
\mu<t'-t<\lambda,~[t,t']\subset I,~\forall s\in [t,t']~|x(s)-x(t)|<\lambda,~
\left|\frac{x(t')-x(t)}{t'-t}-F(t,x(t))\right|<\lambda.$$
By overspill there exists also $\lambda_0\in A_\mu$ with
$0<\lambda_0\simeq 0$. Thus, for all $t\in[0,L]$, there is
$t'\in I$ such that $\mathcal{P}_\mu(t,t',\lambda_0)$ holds.
Applying now the axiom of choice to obtain a function
$c:[0,L]\rightarrow I$ such that $c(t)=t'$, that is,
$\mathcal{P}_\mu(t,c(t),\lambda_0)$ holds for all $t\in[0,L]$.
Since \mbox{$c(t)-t>\mu$} for all $t\in[0,L]$, there are a
positive integer $N$ and an infinitesimal partition $\{t_n:
n=0,\ldots,N+1\}$ of $[0,L]$ such that $t_0=0$, $t_N\leq
L<t_{N+1}$ and $t_{n+1}=c(t_n)$. Finally, the conclusion 
follows from Lemma~\ref{lemma32}.
\end{proof}
%%%%%%%%%%%

\subsubsection{Proof of Theorem \ref{theorem32}}

Let $L>0$ be standard in $J$. Fix $\rho>0$ to be standard and let
$W$ be the (standard) tubular neighborhood around
$\Gamma=\{y(t):t\in[0,L]\}$ given by
\begin{equation*}
W=\{z\in\mathbb{R}^d\ /\ \exists t\in[0,L]: |z-y(t)|\leq\rho\}.
\end{equation*}
Let $A$ be the nonempty ($0\in A$) subset of $[0,L]$ defined by
\begin{equation*}
A=\{L_1\in[0,L]\ /\ [0,L_1]\subset I\hbox{ and } \{x(t):
t\in[0,L_1]\}\subset W\}.
\end{equation*}
The set $A$ is  bounded above by $L$. Let $L_0$ be the upper bound
of $A$ and let $L_1\in A$ be such that $L_0-\mu<L_1\leq L_0$.
Since $\{x(t):t\in[0,L_1]\}\subset W$,  the function $x$ is
limited on $[0,L_1]$.

Taking hypothesis (b) into account, we now apply Lemma
\ref{lemma33} to obtain, for  any standard real number $T$ such
that $0<T\leq L_1$,
\begin{equation}
\label{equation36}
 x(t)\simeq y(t),\quad \forall t\in[0,T].
\end{equation}
By overspill approximation  (\ref{equation36})  still holds for
some $T\simeq L_1$. Next, by Lemma \ref{lemma31} and the
continuity of $y$ we have
\begin{equation*}
 x(t)\simeq x(T) \quad \hbox{and} \quad
y(t)\simeq y(T),\quad \forall t\in[T,L_1].
\end{equation*}
Combining this with (\ref{equation36}) yields
\begin{equation}
\label{equation37}
 x(t)\simeq y(t),\quad \forall t\in[0,L_1].
\end{equation}
Moreover, by hypothesis (a) there exists $L_1'>L_1+\mu$ such that
$[L_1,L_1']\subset I$ and
\begin{equation}
\label{equation38} x(t)\simeq \ y(t),\quad \forall t\in[L_1,L_1'].
\end{equation}
Together (\ref{equation37}) and (\ref{equation38}) show  that
$x(t)\simeq y(t)$ for all $t\in[0,L_1']$.

It remains to verify that $L\leq L_1'$. If this is not true, then
$[0,L_1']\subset I$ and $\{x(t): t\in[0,L_1']\}\subset W$ imply
$L_1'\in A$. This contradicts the fact that $L_1'>L_0$. So the
proof is complete.
%%%%%%%%%

\subsection{The Stroboscopic Method for RFDEs}

Let $\phi\in\mathcal{C}$ be standard and let $F:\mathbb{R}_+\times
\mathcal{C}\rightarrow\mathbb{R}^d$ be a standard and continuous
function. Let $I$ be some subset of $\mathbb{R}$ and let
$x:I\rightarrow\mathbb{R}^d$ be a function such that
$[-r,0]\subset I$ and $x_0\simeq\phi$.

%%%%%%%%%%%%%%%%%%
\begin{definition}[\emph{F}-Stroboscopic property]
\label{definition2} The function $x$ is said to satisfy the
\emph{F}-Stroboscopic property if there exists $\mu>0$ such that,
for all limited $t\geq 0$ in $I$ satisfying $[0,t]\subset I$ and
$x(s)$ and $F(s,x_s)$ are limited for all $s\in[0,t]$, there
exists $t'\in I$ such that $\mu<t'-t\simeq 0$, $[t,t']\subset I$,\
$x(s)\simeq x(t)$ for all $s\in [t,t']$ and $\displaystyle
\frac{x(t')-x(t)}{t'-t}\simeq F(t,x_t). $
\end{definition}
%%%%%%%%%%%%%%%%

In the same manner as in Section \ref{Notations and Main Results},
for $r=0$ we identify the Banach space $\mathcal{C}$ with
$\mathbb{R}^d$ and $x_t$ with $x(t)$. By continuity property of $F$, 
we have then $x(s)$ is
limited for all $s\in[0,t]$ implies that $F(s,x(s))$ is limited
for all $s\in[0,t]$. So, Definition~\ref{definition1} is a
particular case of Definition \ref{definition2}. In the following result
we assert that a function which satisfies the \emph{F}-stroboscopic
property can be approximated by a solution of the RFDE
\begin{equation}
\label{equation39}
 \dot{y}(t) =F(t,y_t),\qquad y_0=\phi.
\end{equation}

%%%%%%%%%%%%%%%
\begin{theorem}[Stroboscopic Lemma for RFDEs]
\label{theorem33} Suppose that\\
{\rm (a)}
The function $x$ satisfies the F-stroboscopic property.\\
{\rm (b)}
The initial value problem (\ref{equation39}) has a unique solution
$y$. Let $J=[-r,\omega)$, $0<\omega\leq\infty$, be its  maximal
interval of definition.\\
Then, for every standard and positive $L\in J$, $[-r,L]\subset I$
and the approximation $x(t)\simeq y(t)$ holds for all $t\in
[-r,L]$.
\end{theorem}
%%%%%%%%%%%%%

To prove Stroboscopic Lemma for RFDEs we need first to establish
the following preliminary  lemmas.

\subsubsection{Preliminaries}

%%%%%%%%%%%%%
\begin{lemma}
\label{lemma34} Let $L>0$ be limited such that $[0,L]\subset I$.
Suppose that\\
{\rm(i)}
$x(t)$ and $F(t,x_t)$ are limited for all $t\in[0,L]$.\\
{\rm(ii)}
There exist some positive integer $N$ and some infinitesimal
partition $\{t_n:n=0,\ldots,N+1\}$ of $[0,L]$ such that $t_0=0$,
$t_N\leq L<t_{N+1}$ and, for $n=0,\ldots,N$, $t_{n+1}\simeq t_n$,
$x(t)\simeq x(t_n)$ for all $t\in [t_n,t_{n+1}]$, and
$\displaystyle \frac{x(t_{n+1})-x(t_n)}{t_{n+1}-t_n}\simeq F(t_n,x_{t_n})$.\\
Then the function $x$ is S-continuous at each point in $[0,L]$.
\end{lemma}
%%%%%%%%%%%
\begin{proof}
Let $t,t'\in[0,L]$ with  $t\leq t'$ and $t\simeq t'$. Let
$p,q\in\{0,\ldots,N\}$ be such that $t\in[t_p,t_{p+1}]$ and
$t'\in[t_q,t_{q+1}]$ with $t_p\simeq t_q$. Define
\begin{equation*}
\eta_n  = \frac{x(t_{n+1})-x(t_n)}{t_{n+1}-t_n}-
F(t_n,x_{t_n}),\quad \forall n\in\{p,\ldots,q-1\}.
\end{equation*}
If $\eta$ and $m$ are the respective maximum values of $|\eta_n|$
and $|F(t_n,x_{t_n})|$ for  $n=p,\ldots,q-1$, then we have
\begin{equation}
\label{equation310}
x(t_q)-x(t_p)  = 
\displaystyle\sum_{n=p}^{q-1}(x(t_{n+1})-x(t_n))
    = \displaystyle \sum_{n=p}^{q-1}
(t_{n+1}-t_n)[F(t_n,x_{t_n})+\eta_n]
\end{equation}
Now, $\eta_n\simeq 0$ for all $n\in\{p,\ldots,q-1\}$ implies
$\eta\simeq 0$. On the other hand, in view of hypothesis (i),
$m=|F(t_s,x_{t_s})|$ for some $s\in\{p,\dots,q-1\}$, is limited.
Then (\ref{equation310}) yields
\begin{equation*}
 |x(t')-x(t)|\simeq |x(t_q)-x(t_p)| \leq (m+\eta)(t_q-t_p) \simeq 0
\end{equation*}
 which shows the S-continuity of $x$ at each point in $[0,L]$.
The proof is complete.
\end{proof}
%%%%%%%%%%%

If the real number $L$ in Lemma \ref{lemma34} is standard one
obtains more precise information about the function $x$.

%%%%%%%%%%%%%
\begin{lemma}
\label{lemma35} Let $L>0$ be standard such that $[0,L]\subset I$.
Suppose that conditions (i) and (ii) in Lemma \ref{lemma34} are
satisfied. Then the standard function
$y:[-r,L]\rightarrow\mathbb{R}^d$ defined by
\begin{equation*}
y_0=\phi,\qquad\mbox{and}\qquad y(t)={}^o\big(x(t)\big)\mbox{ for all standard }t\in[0,L],
\end{equation*}
 is a solution of (\ref{equation39}) and
satisfies
\begin{equation}
\label{equation311}
 x(t)\simeq y(t),\quad \forall t\in[0,L].
\end{equation}
\end{lemma}
%%%%%%%%%%%
\begin{proof}
To prove the lemma one first note that the continuity of $y$  on
$[0,L]$ follows by the same argument as in Lemma \ref{lemma32}. As
a consequence we get the approximation (\ref{equation311}) and
then, for all $t\in[0,L]$, $x_t$ is nearstandard and $x_t\simeq
y_t$. Now it remains to show that $y$ satisfies the integral
equation
\begin{equation*}
y(t)=\phi(0)+\int_0^t F(s,y_s)ds, \quad \forall t\in[0,L].
\end{equation*}
Since the proof does not differ from this one in Step 2 of the
proof of Lemma~\ref{lemma32}, it is omitted.
\end{proof}
%%%%%%%%%%%

From Lemma \ref{lemma35} we deduce the result below.

%%%%%%%%%%%%%
\begin{lemma}
\label{lemma36} Let $L>0$ be standard  such that $[0,L]\subset I$.
Suppose that\\
{\rm(i)}
$x(t)$ and $F(t,x_t)$ are limited for all $t\in[0,L]$.\\
{\rm(ii)}
There exists $\mu>0$ such that, for all $t\in [0,L]$, there exists
$t'\in I$ such that $\mu<t'-t\simeq 0$, $[t,t']\subset I$,
$x(s)\simeq x(t)$ for all $s\in [t,t']$, and 
$\displaystyle\frac{x(t')-x(t)}{t'-t}\simeq F(t,x_t)$.\\
Then the  function $x$ is S-continuous on $[0,L]$ and its shadow
is a solution of (\ref{equation39}) and satisfies approximation
(\ref{equation311}).
\end{lemma}
%%%%%%%%%%%
\begin{proof}
As in the proof of Lemma \ref{lemma33}, we obtain a function
$c:[0,L]\rightarrow I$ satisfying, for all $t\in[0,L]$,
$$\mu<c(t)-t\simeq 0,~[t,c(t)]\subset I,~\forall s\in [t,c(t)]~x(s)\simeq x(t),~
\displaystyle\frac{x(c(t))-x(t)}{c(t)-t}\simeq F(t,x_t).$$
If we let $t_0=0$ and $t_{n+1}=c(t_n)$ for $n=0,\ldots,N$, where
the integer $N$ is such that $t_N\leq L<t_{N+1}$, the conclusion
follows by applying Lemma \ref{lemma35}.
\end{proof}
%%%%%%%%%%%

\subsubsection{Proof of Theorem \ref{theorem33}}
Let $L>0$ be standard in $J$ and $\Gamma=\{y_t:t\in[0,L]\}$. Since
$F$ is standard and continuous, and $[0,L]\times\Gamma$ is a
standard compact subset of $\mathbb{R}^d\times\mathcal{C}$, there
exists $\rho>0$ and standard such that $F$ is limited on
$[0,L]\times W$, where $W$ is the standard tubular neighborhood
around $\Gamma$ defined by
\begin{equation*}
W=\{z\in\mathcal{C}\ /\ \exists t\in[0,L]: |z-y_t|\leq\rho\}.
\end{equation*}

Now, since the set
\begin{equation*}
A=\{L_1\in[0,L]\ /\ [0,L_1]\subset I\hbox{ and } \{x_t:
t\in[0,L_1]\}\subset W\}
\end{equation*}
 is nonempty ($0\in A$) and bounded
above  by $L$, there exists $L_1\in A$ such that $L_0-\mu<L_1\leq
L_0$, where $L_0=\sup A$. However $L_1\in A$ implies
\begin{equation*}
[0,L_1]\times\{x_t:t\in[0,L_1]\}\subset[0,L]\times W
\end{equation*}
 and then
$x(t)$ and $F(t,x_t)$ are limited for all $t\in[0,L_1]$.

If $T$ is standard and  $0<T\leq L_1$, according to Lemma
\ref{lemma36} the shadow on $[0,T]$ of the function $x$ is a
solution of (\ref{equation39}). In view of hypothesis (b), this
shadow coincides with $y$ on $[-r,T]$. Also, we have
\begin{equation*}
 x(t)\simeq y(t),\quad \forall t\in[0,T].
\end{equation*}
By overspill the property above holds for some $T\simeq L_1$. On
the other hand, due to the S-continuity of $x$ at each point in
$([T,L_1]\subset)[0,L_1]$ (see Lemma~\ref{lemma34}) and the
continuity of $y$, we have
\begin{equation*}
x(t)\simeq x(T) \quad \hbox{and} \quad y(t)\simeq y(T),\quad
\forall t\in[T,L_1]
\end{equation*}
which achieves to prove that
\begin{equation}
\label{equation312}
 x(t)\simeq y(t),\quad \forall t\in[0,L_1].
\end{equation}
By hypothesis (a)
\begin{equation}
\label{equation313} x(t)\simeq \ y(t),\quad \forall t\in[L_1,L_1']
\end{equation}
for some $L_1'$ such that $L_1'>L_1+\mu$ and $[L_1,L_1']\subset
I$.
Taking into account that $x_0\simeq\phi=y_0$ and combining
(\ref{equation312}) and (\ref{equation313}), we conclude that
$x_t\simeq y_t$ for all $t\in[0,L_1']$. Assume that $L_1'\leq L$.
Then $[0,L_1']\subset I$ and $\{x_t: t\in[0,L_1']\}\subset W$
imply $L_1'\in A$, which is absurd since $L_1'>L_0$. Thus
$L_1'>L$.
Finally, for any standard $L\in J$ we have shown that $x(t)\simeq
y(t)$ for all $t\in[0,L]\subset[0,L_1']$.
This completes the proof of the theorem.
%%%%%%%

\section{Proofs of the Results}
\label{Proofs of the Results}

We prove Theorems \ref{theorem21}, \ref{theorem24} and
\ref{theorem28} within IST. By {\rm transfer} it suffices to prove
those results for {\em standard data} $f$, $x_0$ and $\phi$. We
will do this by applying Stroboscopic Lemma for ODEs (Theorem
\ref{theorem32}) in both cases of Theorems \ref{theorem21} and
\ref{theorem24}, and Stroboscopic Lemma for RFDEs
(Theorem~\ref{theorem33}) in case of Theorem \ref{theorem28}. For
this purpose we need first to translate all conditions (C1) and
(C2) in Section \ref{Averaging for ODEs}, and (H1), (H2) and (H3)
in Section~\ref{Averaging for RFDEs in normal form} into their
external forms and then prove some technical lemmas.
In the external formulas, we  use the following abbreviations \cite{Nelson}~:
$$\forall^{\rm st}A\mbox{ for }\forall x({\rm st}x\Rightarrow A)
\qquad\mbox{ and }\qquad
\exists^{\rm st}A\mbox{ for }\exists x({\rm st}x~\& A).$$
Let $f:\mathbb{R}_+\times\mathcal{C}\rightarrow\mathbb{R}^d$ be a
standard and continuous function, where
$\mathcal{C}=\mathcal{C}([-r,0],\mathbb{R}^d)$ and $r\geq 0$. We
recall that when $r=0$, $\mathcal{C}$ is identified with
$\mathbb{R}^d$. The external formulations of conditions (C1) and
(C2) are:
%%%%%%%%%%%%%%%%%%%%%%%%%%%%%%%%%%%%%%%%%%%%%%%%%%%
\begin{enumerate}
    \item[(C1')]
$\forall^{\rm st} x\in\mathbb{R}^d$ $\forall x'\in\mathbb{R}^d$
$\forall t\in\mathbb{R}_+$ $\big(x'\simeq x \Rightarrow
f(t,x')\simeq f(t,x)\big).$
    \item[(C2')]
$\exists^{\rm st}F:\mathbb{R}^d\rightarrow\mathbb{R}^d$ $ \forall^{\rm st}
x\in\mathbb{R}^d$ $\forall R\simeq +\infty$
$F(x)\simeq\displaystyle\frac{1}{R}\int_0^Rf(t,x)dt$.
\end{enumerate}
The external formulation of conditions (H1), (H2) and (H3) are,
respectively:
\begin{enumerate}
    \item[(H1')]
$\forall^{\rm st} x\in\mathcal{C}$ $\forall x'\in\mathcal{C}$ $\forall
t\in\mathbb{R}_+$ $\big(x'\simeq x \Rightarrow f(t,x')\simeq
f(t,x)\big).$
    \item[(H2')]
$\forall x\in\mathcal{C}$ and limited $\forall t\in\mathbb{R}_+$,
$f(t,x)$ is limited in $\mathbb{R}^d$.
    \item[(H3')]
$\exists^{\rm st}F:\mathcal{C}\rightarrow\mathbb{R}^d$ $ \forall^{\rm st}
x\in\mathcal{C}$ $\forall R\simeq +\infty$
$F(x)\simeq\displaystyle\frac{1}{R}\int_0^Rf(t,x)dt$.
\end{enumerate}

\subsection{Technical Lemmas}

In Lemmas \ref{lemma41} and \ref{lemma42} below we formulate some
properties of the average of the function $f$ (i.e. the function
$F$ defined in (C2) and (H3)).

%%%%%%%%%%%%%
\begin{lemma}
\label{lemma41} Suppose that the function $f$ satisfies conditions
(C1) and (C2) when $r=0$ and conditions (H1) and (H3) when $r>0$.
Then the function $F$ in (C2) or (H3) is continuous and satisfies
\begin{equation*}
F(x)\simeq\frac{1}{R}\int_0^R f(t,x)dt
\end{equation*}
 for all nearstandard $x\in\mathcal{C}$ and all $R\simeq +\infty$.
\end{lemma}
%%%%%%%%%%%
\begin{proof} The proof is the same in both cases $r=0$ and $r>0$. So,
there is no restriction to suppose that $r=0$. Let $x,
{}^ox\in\mathbb{R}^d$ be such that ${}^ox$ is standard and
${}^ox\simeq x$. Fix $\delta>0$ to be infinitesimal. By condition
(C2)
\begin{equation*}
\left|F(x)-\frac{1}{T}\int_0^Tf(t,x)dt\right|<\delta,\quad \forall
T>T_0
\end{equation*}
 for some $T_0>0$.  Hence there exists $T\simeq+\infty$ such
that
\begin{equation*}
 F(x)\simeq\frac{1}{T}\int_0^T f(t,x)dt.
\end{equation*}
By condition (C1') we have $f(t,x)\simeq f(t,{}^ox)$ for all
$t\in\mathbb{R}_+$. Therefore
\begin{equation*}
 F(x)\simeq\frac{1}{T}\int_0^T f(t,{}^ox)dt.
\end{equation*}
By condition (C2') we deduce that $F(x)\simeq F({}^ox)$. Thus $F$
is continuous. Moreover, for all $T\simeq+\infty$, we have
\begin{equation*}
 F(x)\simeq F({}^ox)\simeq\frac{1}{T}\int_0^T f(t,{}^ox)dt\simeq\frac{1}{T}\int_0^T f(t,x)dt.
\end{equation*}
So, the proof is complete.
\end{proof}
%%%%%%%%%%%
\begin{lemma}
\label{lemma42}  Suppose that the function $f$ satisfies
conditions (C1) and (C2) when $r=0$ and conditions (H1) and (H3)
when $r>0$. Let $F$ be as in (C2) or (H3). Let $\varepsilon>0$ be
infinitesimal. Then, for all limited $t\in\mathbb{R}_+$ and all
nearstandard $x\in\mathcal{C}$, there exists
$\alpha=\alpha(\varepsilon,t,x)$ such that $0<\alpha\simeq 0$,
$\varepsilon/\alpha\simeq 0$ and
\begin{equation*}
\frac{\varepsilon}{\alpha}\int_{t/\varepsilon}^{t/\varepsilon+T\alpha/\varepsilon}
f(\tau,x)d\tau\simeq TF(x), \quad \forall T\in [0,1].
\end{equation*}
\end{lemma}
%%%%%%%%%%%
\begin{proof}
The proof is the same in both cases $r=0$ and $r>0$. Let $t$ be
limited in $\mathbb{R}_+$ and let $x$ be nearstandard in
$\mathcal{C}$. We denote for short $g(r)= f(r,x)$. Let
$T\in[0,1]$. We consider the following two cases.

Case 1: $t/\varepsilon$ is limited. Let $\alpha>0$ be such that
$\varepsilon/\alpha\simeq 0$. If $T\alpha/\varepsilon$ is limited
then we have $T\simeq0$ and
\begin{equation*}
\frac{\varepsilon}{\alpha}\int_{t/\varepsilon}^{t/\varepsilon+T\varepsilon/\alpha}g(r)dr\simeq0\simeq
TF(x).
\end{equation*}
 If $T\alpha/\varepsilon\simeq+\infty$ we write
\begin{equation*}
\displaystyle
\frac{\varepsilon}{\alpha}\int_{t/\varepsilon}^{t/\varepsilon+T\alpha/\varepsilon}g(r)dr=
\displaystyle\left(T+\frac{t}{\alpha}\right)
\frac{1}{t/\varepsilon+T\alpha/\varepsilon}\int_0^{t/\varepsilon+T\alpha/\varepsilon}g(r)dr-
\frac{\varepsilon}{\alpha}\int_0^{t/\varepsilon}g(r)dr.
\end{equation*}
By Lemma \ref{lemma41} we have
\begin{equation*}
\frac{1}{t/\varepsilon+T\alpha/\varepsilon}\int_0^{t/\varepsilon+T\alpha/\varepsilon}g(r)dr\simeq
F(x).
\end{equation*}
 Since
$\frac{\varepsilon}{\alpha}\int_0^{t/\varepsilon}g(r)dr\simeq0$
and $t/\alpha\simeq 0$, we have
\begin{equation*}
\frac{\varepsilon}{\alpha}\int_{t/\varepsilon}^{t/\varepsilon+T\alpha/\varepsilon}g(r)dr\simeq
TF(x).
\end{equation*}
This approximation is satisfied for all $\alpha>0$ such that
$\varepsilon/\alpha\simeq 0$. Choosing then $\alpha$ such that
$0<\alpha\simeq 0$ and $\varepsilon/\alpha\simeq 0$ gives the
desired result.

Case 2: $t/\varepsilon$ is unlimited. Let $\alpha>0$. We have
\begin{equation}
\label{equation41}
\displaystyle
\frac{\varepsilon}{\alpha}\int_{t/\varepsilon}^{t/\varepsilon+T\alpha/\varepsilon}g(r)dr
 =  T\eta(\alpha)
+\frac{t}{\alpha}\left[\eta(\alpha)-\eta(0)\right],
\end{equation}
where
\begin{equation*}
\eta(\alpha)=\frac{1}{t/\varepsilon+T\alpha/\varepsilon}\int_0^{t/\varepsilon+T\alpha/\varepsilon}g(r)dr.
\end{equation*}
By Lemma \ref{lemma41} we have $\eta(\alpha)\simeq F(x)$ for all
$\alpha\geq 0$.
Return to (\ref{equation41}) and assume that $\alpha$ is not
infinitesimal. Then
\begin{equation}
\label{equation42}
\displaystyle
\frac{\varepsilon}{\alpha}\int_{t/\varepsilon}^{t/\varepsilon+T\alpha/\varepsilon}g(r)dr\simeq TF(x),
\end{equation}
By overspill (\ref{equation42}) holds  for some
$\alpha\simeq 0$ which can be chosen such that
$\varepsilon/\alpha\simeq 0$.
This completes the proof of the lemma.
\end{proof}
%%%%%%%%%%%

\subsection{Proof of Theorem \ref{theorem21}}
\label{Proof1}

Assume that $x_0$ and $L$ are standard. To prove Theorem
\ref{theorem21} is equivalent to show that, for every
infinitesimal $\varepsilon>0$, every solution $x$ of
(\ref{equation21}) is defined at least on $[0,L]$ and satisfies
$x(t)\simeq y(t)$ for all $t\in [0,L]$.
We need first to prove the following result which discuss some
properties of solutions of a certain ODE needed in the sequel.
%%%%%%%%%%%%%%%%%%%%%%%%%%%%%%%%%%%%%%%%%%%%%%%%%
\begin{lemma}
\label{lemma43} Let $g:\mathbb{R}_+\times
\mathbb{R}^d\rightarrow\mathbb{R}^d$ and
$h:\mathbb{R}_+\rightarrow\mathbb{R}^d$ be continuous functions.
Let $x_0$  be limited in $\mathbb{R}^d$. Suppose that
\begin{enumerate}
    \item[\rm (i) ]
$g(t,x)\simeq h(t)$ holds for all $t\in[0,1]$ and all limited
$x\in\mathbb{R}^d$.
    \item[\rm (ii) ]
$\displaystyle\int_0^t h(s)ds$ is limited for all $t\in[0,1]$.
\end{enumerate}
Then any solution $x$ of the initial value problem
\begin{equation*}
 \dot x=g(t,x),\ t\in[0,1];\qquad x(0)=x_0
\end{equation*}
is defined and limited on $[0,1]$ and satisfies
\begin{equation*}
x(t)\simeq x_0+\int_0^t h(s)ds,\quad \forall t\in[0,1].
\end{equation*}
\end{lemma}
%%%%%%%%%%%
\begin{proof} By overspill there exists $\omega\simeq+\infty$ such
that the approximation in hypothesis (i) holds for all $t\in[0,1]$
and all $x\in B(0,\omega)$, where
\hbox{$B(0,\omega)\subset\mathbb{R}^d$} is the ball of center $0$
and radius $\omega$. Assume  that $x(t')\simeq\infty$ for some
$t'\in[0,1]$. Let $t\in[0,1]$ be such that $t\leq t'$, $x(t)\simeq
\infty$ and $x(s)\in B(0,\omega)$ for all $s\in[0,t]$. Then we
have
\begin{equation*}
 x(t)=x_0+\int_0^t
g(s,x(s))ds\simeq x_0+\int_0^t h(s)ds
\end{equation*}
whence, in view of hypothesis (ii), $x(t)$ is limited; this is a
contradiction. Therefore $x(t)$ is defined and limited for all
$t\in[0,1]$.
\end{proof}
%%%%%%%%%%%

Let us now  prove Theorem \ref{theorem21}. Fix $\varepsilon>0$
to be infinitesimal and let $x:I\rightarrow\mathbb{R}^d$ be a
maximal solution of~(\ref{equation21}). We claim that $x$
satisfies the \emph{F}-stroboscopic property. To see this, let
$t_0\in I$ such that $t_0\geq 0$ and limited, and $x(t)$ is
limited for all $t\in[0,t_0]$. By Lemma~\ref{lemma42} there exists
$\alpha=\alpha(\varepsilon,t_0,x(t_0))$ such that $0<\alpha\simeq
0$, $\varepsilon/\alpha\simeq 0$ and
\begin{equation}
\label{equation43}
\frac{\varepsilon}{\alpha}\int_{t_0/\varepsilon}^{t_0/\varepsilon+T\alpha/\varepsilon}
 f(t,x(t_0))dt\simeq TF(x(t_0)),\quad \forall T\in[0,1].
\end{equation}

Introduce the function
\begin{equation*}
X(T)=\frac{x(t_0+\alpha T)-x(t_0)}{\alpha},\quad T\in[0,1].
\end{equation*}
Differentiating and substituting the above into (\ref{equation21})
gives, for $T\in[0,1]$,
\begin{equation}
\label{equation44}
 \frac{dX}{dT}(T) =f\left(\frac{t_0}{\varepsilon}+
\frac{\alpha}{\varepsilon}T,x(t_0)+\alpha X(T)\right).
\end{equation}
By (C1') and Lemma \ref{lemma43}  the function $X$, as a solution
of (\ref{equation44}),  is defined and limited on $[0,1]$ and, for
$T\in[0,1]$,
\begin{equation*}
 X(T) \simeq
 \int_0^T
f\left(\frac{t_0}{\varepsilon}+\frac{\alpha}{\varepsilon}t,x(t_0)\right)dt
    = \displaystyle   \frac{\varepsilon}{\alpha}
  \int_{t_0/\varepsilon}^{t_0/\varepsilon+T\alpha/\varepsilon} f(t,x(t_0))dt.
\end{equation*}
Using now (\ref{equation43}) this leads to the approximation
\begin{equation*}
X(T)\simeq \displaystyle TF(x(t_0)), \quad\forall T\in[0,1].
\end{equation*}
Define $t_1=t_0+\alpha$ and set $\mu=\varepsilon$.  Then
$\mu<\alpha=t_1-t_0\simeq 0$, $[t_0,t_1]\subset I$ and
$x(t_0+\alpha T)=x(t_0)+\alpha X(T)\simeq x(t_0)$ for all
$T\in[0,1]$, that is, $x(t)\simeq x(t_0)$ for all $t\in[t_0,t_1]$,
whereas
\begin{equation*}
 \frac{x(t_1)-x(t_0)}{t_1-t_0}=X(1)\simeq
F(x(t_0)),
\end{equation*}
 which shows the claim. Finally, by (C3) and Theorem
\ref{theorem32}, the solution $x$ is defined at least on $[0,L]$
and satisfies $x(t)\simeq y(t)$ for all $t\in [0,L]$.
The proof of Theorem \ref{theorem21} is complete.
%%%%%%%

\subsection{Proof of Theorem \ref{theorem24}}

Assume that $\phi$ and $L$ are standard. All we have to prove is
that, when $\varepsilon>0$ is infinitesimal, every solution $x$ of
(\ref{equation23}) is defined at least on $[-\varepsilon r,L]$ and
satisfies $x(t)\simeq y(t)$ for all $t\in [0,L]$.
We start by showing some auxiliary result which is needed in the
proof of Theorem \ref{theorem24}.

%%%%%%%%%%%%%
\begin{lemma}
\label{lemma44} Let
$g:\mathbb{R}_+\times\mathcal{C}\rightarrow\mathbb{R}^d$ be a
continuous function. Suppose that, for all $t\in\mathbb{R}_+$ and
all $x\in\mathcal{C}$, $t$ and $x$ limited imply that $g(t,x)$ is
limited. Let  $\phi\in\mathcal{C}$ be standard. Let
$x:I\rightarrow\mathbb{R}^d$ be a maximal solution of the initial
value problem
\begin{equation}
\label{equation45} \dot{x}(t) = g(t,x_{t,\varepsilon}),\qquad
x(t)=\phi(t/\varepsilon),\ t\in[-\varepsilon r,0].
\end{equation}
Let $t_0\geq 0$ be limited in $I$ such that $x$ is limited on
$[0,t_0]$. Then the solution $x$ is
\begin{enumerate}
    \item[\rm(i) ]
S-continuous at each $t\in[0,t_0]$.
    \item[\rm(ii) ]
defined and  limited at each $t\geq t_0$ such that $t\simeq t_0$.
\end{enumerate}
\end{lemma}
%%%%%%%%%%%
\begin{proof}
(i) Let $t\in[0,t_0]$. We will show that $x$ is S-continuous at
$t$. If $t'\in[0,t_0]$ is such that $t\leq t'$ and $t\simeq t'$,
the integral equation for the solutions of (\ref{equation45})
implies
\begin{equation*}
|x(t')-x(t)|\leq \int_t^{t'}|g(s,x_{s,\varepsilon})|ds \leq
(t'-t)\sup_{s\in[t,t']}|g(s,x_{s,\varepsilon})|.
\end{equation*}
Since $x_{s,\varepsilon}(\theta)=x(s+\varepsilon\theta)$ is
limited for all $\theta\in[-r,0]$ and all $s\in[t,t']$, we have
$x_{s,\varepsilon}$ is limited for all $s\in[t,t']$. Now in view
of assumptions on $x$ and~$g$,
$\sup_{s\in[t,t']}|g(s,x_{s,\varepsilon})|$ is limited so that
$x(t')\simeq x(t)$, which is the desired result.

(ii) Let $I=[-\varepsilon r,b)$, $0<b\leq\infty$. We assume by
contradiction that $x$ is not defined for all $t\geq t_0$ such
that $t\simeq t_0$. Then $b\simeq t_0$. We have
$x(t')\simeq\infty$ for some $t'\in[t_0,b)$. Indeed, taking into
account that $x$ is limited on $[-\varepsilon r,0]$ since
$x([-\varepsilon r,0])=\phi([-r,0])$, if it was also limited on
$[t_0,b)$ then $\lim_{t\rightarrow b}x(t)$ exists and $x$ can be
extended to a continuous function on $[-\varepsilon r,b]$ by
setting $x(b)=\lim_{t\rightarrow b}x(t)$. Consequently,
$x_{b,\varepsilon}\in\mathcal{C}$ and then one can find a solution
of (\ref{equation45}) through the point $(b,x_{b,\varepsilon})$ to
the right of $b$, which contradicts the noncontinuability
hypothesis on $x$. Now, by the continuity of $x$ there exists
$t\in[t_0,b)$ such that $x$ is limited on $[t_0,t]$ and
$x(t)\not\simeq x(t_0)$. Using Part (i) of the lemma,  $x$ is
S-continuous at each point in $[0,t]$ and $x(t)\simeq x(t_0)$
since $t\simeq t_0$. This is a contradiction and proves the first
claim of Part (ii) of the lemma.
If the second claim of Part (ii) is not true, then
$x(t')\simeq\infty$ for some $t'\in[t_0,b)$ with $t'\simeq t_0$.
Again using continuity, $x$ is  limited on $[t_0,t]$ and
$x(t)\not\simeq x(t_0)$ for some $t\in[t_0,b)$. The same reasoning
as above gives a contradiction and proves that $x(t)$ is limited
for all $t\geq t_0$ such that $t\simeq t_0$.
This completes the proof of the lemma.
\end{proof}
%%%%%%%%%%%

The proof of Theorem \ref{theorem24} is as follows. Let
$\varepsilon>0$ be infinitesimal. Let $x$ be a maximal solution of
(\ref{equation23}) defined on $I$, an interval of $\mathbb{R}$.
Let $t_0\in I$ such that $t_0\geq 0$ and limited, and $x(t)$ is
limited for all $t\in[0,t_0]$. Since $x(t_0)$ is limited, it is
nearstandard and so is $\tilde{x}^{t_0}$ where
$\tilde{x}^{t_0}\in\mathcal{C}$ is defined by
$\tilde{x}^{t_0}(\theta)=x(t_0)$ for all $\theta\in[-r,0]$. Now we
apply Lemma \ref{lemma42} to obtain some constant
$\alpha=\alpha(\varepsilon,t_0,\tilde{x}^{t_0})$ such that
$0<\alpha\simeq 0$, $\varepsilon/\alpha\simeq 0$ and
\begin{equation}
\label{equation46}
\frac{\varepsilon}{\alpha}\int_{t_0/\varepsilon}^{t_0/\varepsilon+T\alpha/\varepsilon}
 f(t,\tilde{x}^{t_0})dt\simeq TF(\tilde{x}^{t_0})=TG(x(t_0)),\quad \forall
 T\in[0,1].
\end{equation}

By Lemma \ref{lemma44} the solution $x$ is defined for all $t\geq
t_0$ and $t\simeq t_0$. Hence one can consider the function
\begin{equation*}
 X(\theta,T)=\frac{x(t_0+\alpha
T+\varepsilon\theta)-x(t_0)}{\alpha},\quad \theta\in[-r,0], \
T\in[0,1].
\end{equation*}
We have, for $T\in[0,1]$,
\begin{equation*}
X(0,T)=\frac{x(t_0+\alpha T)-x(t_0)}{\alpha},\quad x_{t_0+\alpha
T,\varepsilon}=\tilde{x}^{t_0}+\alpha X(\cdot,T)
\end{equation*}
 and therefore
\begin{equation*}
 \frac{\partial X}{\partial T}(0,T)
=f\left(\frac{t_0}{\varepsilon}+
\frac{\alpha}{\varepsilon}T,\tilde{x}^{t_0}+\alpha
X(\cdot,T)\right).
\end{equation*}
 Integration between $0$ and $T$, for
$T\in[0,1]$, yields
\begin{equation}
\label{equation47}
X(0,T)=\int_0^Tf\left(\frac{t_0}{\varepsilon}+\frac{\alpha}{\varepsilon}t,
         \tilde{x}^{t_0}+\alpha X(\cdot,t))\right)dt.
\end{equation}

We now consider the  following two cases:

Case 1: $T\in[0,\varepsilon r/\alpha]$. Note that
$\tilde{x}^{t_0}+\alpha X(\cdot,T)$ is limited for all $T\in[0,1]$
(in particular for all $T\in[0,\varepsilon r/\alpha]$). Using
(H3') and taking into account that $\varepsilon r/\alpha\simeq 0$,
(\ref{equation47}) leads to the approximation
 \begin{equation}
\label{equation48}
 X(0,T)\simeq 0.
 \end{equation}

Case 2: $T\in[\varepsilon r/\alpha,1]$.  By Lemma \ref{lemma44}
the solution  $x$ is S-continuous at each point in
$[0,t_0+\alpha]$ so that, for $\theta\in[-r,0]$,
\begin{equation*}
 \alpha X(\theta,T)=x(t_0+\alpha T+\varepsilon\theta)-x(t_0) \simeq
 0,
\end{equation*}
since  $t_0+\alpha
T+\varepsilon\theta\in[t_0,t_0+\alpha]\subset[0,t_0+\alpha]$ and
$t_0+\alpha T+\varepsilon\theta\simeq t_0$.

Return now to (\ref{equation47}). For $T\in[0,1]$, we write
\begin{equation*}
X(0,T)  =  \left(\int_0^{\varepsilon r/\alpha}+\int_{\varepsilon r/\alpha}^T\right)
f\left(\frac{t_0}{\varepsilon}+\frac{\alpha}{\varepsilon}t,
         \tilde{x}^{t_0}+\alpha X(\cdot,t)\right)dt.
\end{equation*}
Note that $x(t_0)$ is nearstandard implies that, for all
\mbox{$T\in[0,1]$}, $\tilde{x}^{t_0}+\alpha X(\cdot,T)$ is
nearstandard. Using (\ref{equation48}), (H1') and (H2'), we thus
get, for \hbox{$T\in[0,1]$},
\begin{equation*}
\begin{matrix}
X(0,T)  &\simeq  & \displaystyle \int_{\varepsilon
r/\alpha}^Tf\left(\frac{t_0}{\varepsilon}+\frac{\alpha}{\varepsilon}t,
         \tilde{x}^{t_0}\right)dt%\hfill\\[5mm]
\simeq \displaystyle\int_0^T
f\left(\frac{t_0}{\varepsilon}+\frac{\alpha}{\varepsilon}t,
         \tilde{x}^{t_0}\right)dt\hfill\\[5mm]
   & = & \displaystyle  \frac{\varepsilon}{\alpha}
  \int_{t_0/\varepsilon}^{t_0/\varepsilon+T\alpha/\varepsilon} f(t, \tilde{x}^{t_0})dt\hfill
\end{matrix}
\end{equation*}
whence, in view of  (\ref{equation46}),
\begin{equation*}
X(0,T)\simeq \displaystyle TG(x(t_0)).
\end{equation*}
 Defining $t_1=t_0+\alpha$
and setting  $\mu=\varepsilon$, the following properties are true:
$\mu<\alpha=t_1-t_0\simeq 0$, $[t_0,t_1]\subset I$, $x(t_0+\alpha
T)=x(t_0)+\alpha X(0,T)\simeq x(t_0)$ for all $T\in[0,1]$, that
is, $x(t)\simeq x(t_0)$ for all $t\in[t_0,t_1]$, and
\begin{equation*}
 \frac{x(t_1)-x(t_0)}{t_1-t_0}=X(0,1)\simeq
G(x(t_0)).
\end{equation*}
 This proves that $x$ satisfies the
\emph{F}-stroboscopic property. Taking (H4) into account, we
finally apply  Theorem \ref{theorem32} (Stroboscopic Lemma for
ODEs) to obtain the desired result, that is, the solution $x$ is
defined at least on $[-\varepsilon r,L]$ and satisfies $x(t)\simeq
y(t)$ for all $t\in [0,L]$.
The theorem is proved.
%%%%%%%%

\subsection{Proof of Theorem \ref{theorem28}}

Let  $\phi$  and $L$ be standard. To prove Theorem
\ref{theorem28} is equivalent to show that for every infinitesimal
$\varepsilon>0$, every solution $x$ of (\ref{equation25}) is
defined at least on $[-r,L]$ and $x(t)\simeq y(t)$ holds for all
$t\in [0,L]$.
Before this, we first prove the following result.

%%%%%%%%%%%%%
\begin{lemma}
\label{lemma45} Let
$g:\mathbb{R}_+\times\mathcal{C}\rightarrow\mathbb{R}^d$ be a
continuous function. Suppose that, for all $t\in\mathbb{R}_+$ and
all $x\in\mathcal{C}$, $t$ and $x$ limited imply that $g(t,x)$ is
limited. Let  $\phi\in\mathcal{C}$ be standard. Let
$x:I\rightarrow\mathbb{R}^d$ be a maximal solution of the initial
value problem
\begin{equation*}
 \dot{x}(t) =  g(t,x_t),\qquad x_0=\phi.
\end{equation*}
Let $t_0\geq 0$ be limited in $I$ such that $x$ is limited on
$[0,t_0]$. Then
\begin{enumerate}
    \item[\rm(i) ]
$x$ is S-continuous at each $t\in[-r,t_0]$ and $x_t$ is
nearstandard for all $t\in [0,t_0]$.
    \item[\rm(ii) ]
$x$ is defined and  limited at each $t\simeq t_0$.
\end{enumerate}
\end{lemma}
%%%%%%%%%%%
\begin{proof}
To prove (i) we first note  that $x$ is S-continuous at each
\mbox{$t\in[-r,0]$}, since it coincides with the standard and
continuous function $\phi$ on the (standard) interval $[-r,0]$.
Now consider the interval $[0,t_0]$. Let $t\in[0,t_0]$. If
$t'\in[0,t_0]$ is such that $t\leq t'$ and $t\simeq t'$, then
\begin{equation*}
|x(t')-x(t)|\leq \int_t^{t'}|g(s,x_s)|ds \leq
(t'-t)\sup_{s\in[t,t']}|g(s,x_s)|.
\end{equation*}
In view of assumptions on $x$ and $g$, the quantity
$\sup_{s\in[t,t']}|g(s,x_s)|$ is limited so that  $x(t')\simeq
x(t)$. This shows the S-continuity of $x$ at $t$.

It remains to prove that $x_t$ is nearstandard for all $t\in
[0,t_0]$. We have, $x$ is limited and S-continuous at each
$t\in[-r,t_0]$ implies that $x_t$ is limited and S-continuous for
all $t\in [0,t_0]$. So, the desired result follows from
Theorem~\ref{theorem31}.

To prove (ii) we let $I=[-r,b)$, $0<b\leq\infty$, and suppose that
$x$ is not defined for all $t\simeq t_0$, that is, $b\simeq t_0$.
Then there exists $t'\in[t_0,b)$ such that $x(t')\simeq\infty$.
Otherwise, $\lim_{t\rightarrow b}x(t)$ exists and $x$ can be
continued through the point $(b,x_b)$ to the right of $b$, which
contradicts the noncontinuability hypothesis on $x$. Now, by the
continuity of $x$ there exists $t\in[t_0,b)$ such that $x$ is
limited on $[t_0,t]$ and $x(t)\not\simeq x(t_0)$. On the other
hand, by Part~(i) of  the lemma $x$ is S-continuous at each point
in $[-r,t]$. Since $t\simeq t_0$, it follows that $x(t)\simeq
x(t_0)$, which is absurd. This proves that $x(t)$ is defined for
all $t\simeq t_0$.

Suppose now that $x(t)$ is not limited for all $t\simeq t_0$, that
is, $x(t')\simeq\infty$ for some $t'\in[t_0,b)$ with $t'\simeq
t_0$. Again, by the continuity of $x$ there exists
\mbox{$t\in[t_0,b)$} such that $x$ is  limited on $[t_0,t]$ and
$x(t)\not\simeq x(t_0)$. The same argument as above leads to a
contradiction. This proves that $x(t)$ is limited for all $t\simeq
t_0$. Lemma \ref{lemma45} is proved.
\end{proof}
%%%%%%%%%%%

For the proof of Theorem \ref{theorem28} we fix $\varepsilon>0$ to
be infinitesimal and we let $x:I\rightarrow\mathbb{R}^d$ to be a
maximal solution of (\ref{equation25}). We will first show that
$x$ satisfies the \emph{F}-stroboscopic property. Let $t_0\in I$
such that $t_0\geq 0$ and limited, and $x(t)$ and $F(x_t)$ are
limited for all $t\in[0,t_0]$. According to (H2') Lemma
\ref{lemma45} applies. Thus $x_t$ is nearstandard for all
$t\in[0,t_0]$.

Now, applied to $t_0$ and $x_{t_0}$, Lemma \ref{lemma42} gives
\begin{equation}
\label{equation49}
\frac{\varepsilon}{\alpha}\int_{t_0/\varepsilon}^{t_0/\varepsilon+T\alpha/\varepsilon}
 f(t,x_{t_0})dt\simeq TF(x_{t_0}),\quad \forall T\in[0,1]
\end{equation}
for some $\alpha=\alpha(\varepsilon,t_0,x_{t_0})$ such that
$0<\alpha\simeq 0$ and $\varepsilon/\alpha\simeq 0$.

Let $X:[-r,0]\times[0,1]\rightarrow\mathbb{R}^d$ be the function
given by
\begin{equation*}
X(\theta,T)=\frac{x(t_0+\alpha
T+\theta)-x(t_0+\theta)}{\alpha},\quad \theta\in[-r,0], \
T\in[0,1].
\end{equation*}
By Lemma \ref{lemma45} the function $X$ is well defined. It
satisfies, for $T\in[0,1]$,
\begin{equation*}
X(0,T)=\frac{x(t_0+\alpha T)-x(t_0)}{\alpha},\quad x_{t_0+\alpha
T}=x_{t_0}+\alpha X(\cdot,T).
\end{equation*}
Hence, for $T\in[0,1]$,
\begin{equation*}
 \frac{\partial X}{\partial T}(0,T)
=f\left(\frac{t_0}{\varepsilon}+
\frac{\alpha}{\varepsilon}T,x_{t_0}+\alpha X(\cdot,T)\right).
\end{equation*}
Solving this equation gives, for $T\in[0,1]$,
\begin{equation}
\label{equation410}
X(0,T)=\int_0^Tf\left(\frac{t_0}{\varepsilon}+\frac{\alpha}{\varepsilon}t,
         x_{t_0}+\alpha X(\cdot,t)\right)dt.
\end{equation}

Now according to Lemma \ref{lemma45} the solution $x$ is
S-continuous at each point in $[-r,t_0+\alpha]$. Therefore, for
$\theta\in[-r,0]$ and $T\in[0,1]$,  $X(\theta,T)$ satisfies, since
$t_0+\alpha T+\theta\simeq t_0+\theta$,
\begin{equation}
\label{equation411}
 \alpha X(\theta,T)=x(t_0+\alpha T+\theta)-x(t_0+\theta) \simeq 0.
\end{equation}

By (H1') and taking (\ref{equation411}) into account,
(\ref{equation410}) leads, for $T\in[0,1]$, to the approximation
\begin{equation*}
 X(0,T) \simeq \displaystyle
\int_0^T
f\left(\frac{t_0}{\varepsilon}+\frac{\alpha}{\varepsilon}t,x_{t_0}\right)dt
    = \frac{\varepsilon}{\alpha}
  \int_{t_0/\varepsilon}^{t_0/\varepsilon+T\alpha/\varepsilon} f(t,x_{t_0})dt.
\end{equation*}
From (\ref{equation49}) we have
\begin{equation*}
X(0,T)\simeq \displaystyle TF(x_{t_0}),\quad\forall T\in[0,1].
\end{equation*}
Let $t_1=t_0+\alpha$ and set $\mu=\varepsilon$. The instant $t_1$
and the constant $\mu$ are such that: $\mu<\alpha=t_1-t_0\simeq
0$, $[t_0,t_1]\subset I$, $x(t_0+\alpha T)=x(t_0)+\alpha
X(0,T)\simeq x(t_0)$ for all $T\in[0,1]$, that is, $x(t)\simeq
x(t_0)$ for all $t\in[t_0,t_1]$ and
\begin{equation*}
\frac{x(t_1)-x(t_0)}{t_1-t_0}=X(0,1)\simeq F(x_{t_0}),
\end{equation*}
 which form the \emph{F}-stroboscopic property. Finally, using (H6) we
get, by means of Theorem \ref{theorem33} (Stroboscopic Lemma for
RFDEs), the solution $x$ is defined at least on $[-r,L]$ and
satisfies \hbox{$x(t)\simeq y(t)$} for all $t\in [0,L]$. So the
proof is complete.

\end{document}